\documentclass[11pt,twoside,a4paper]{article}
\usepackage{amsmath,amssymb,amsthm}

\usepackage{graphicx,psfrag}
\newtheorem{thm}{Theorem}[section]
\newtheorem{lem}[thm]{Lemma}
\newtheorem{pro}[thm]{Proposition}
\newtheorem{cor}[thm]{Corollary}

\theoremstyle{definition}

\theoremstyle{remark}
\newtheorem{rem}[thm]{Remark}

\newcommand{\R}{\mathbb{R}}

\newcommand{\N}{\mathbb{N}}

\newcommand{\cN}{\mathcal{N}}

\newcommand{\cP}{\mathcal{P}}

\newcommand{\cU}{\mathcal{U}}
\newcommand{\cV}{\mathcal{V}}
\newcommand{\cW}{\mathcal{W}}

\newcommand{\al}{\alpha}

\newcommand{\ga}{\gamma}

\newcommand{\de}{\delta}
\newcommand{\De}{\Delta}

\newcommand{\si}{\sigma}

\newcommand{\la}{\lambda}

\renewcommand{\phi}{\varphi}

\newcommand{\dist}{\operatorname{dist}}
\newcommand{\diam}{\operatorname{diam}}

\newcommand{\hyp}{\operatorname{H}}

\newcommand{\Lip}{\operatorname{Lip}}

\newcommand{\const}{\operatorname{const}}
\newcommand{\Int}{\operatorname{Int}}

\newcommand{\asdim}{\operatorname{asdim}}
\newcommand{\ba}{\operatorname{ba}}
\newcommand{\st}{\operatorname{st}}
\newcommand{\pt}{\operatorname{pt}}

\newcommand{\cone}{\operatorname{Co}}

\newcommand{\mesh}{\operatorname{mesh}}
\newcommand{\cdim}{\operatorname{cdim}}
\newcommand{\andim}{\operatorname{ANdim}}
\newcommand{\cp}{\operatorname{cap}}
\newcommand{\loc}{\operatorname{loc}}
\newcommand{\cploc}{\operatorname{cap}_{\loc}}

\newcommand{\es}{\emptyset}
\renewcommand{\d}{\partial}
\newcommand{\di}{\d_{\infty}}
\newcommand{\set}[2]{\{#1:\,\text{#2}\}}
\newcommand{\sm}{\setminus}
\newcommand{\sub}{\subset}
\newcommand{\sups}{\supset}

\newcommand{\ov}{\overline}

\newcommand{\wh}{\widehat}

\begin{document}

\title{Asymptotic dimension of a hyperbolic space and
capacity dimension of its boundary at infinity}
\author{Sergei Buyalo\footnote{Supported by RFFI Grant
05-01-00939 and Grant NSH-1914.2003.1}}

\date{}
\maketitle

\begin{abstract} We introduce a quasi-symmetry invariant
of a metric space
$Z$
called the capacity dimension,
$\cdim Z$.
Our main result says that for a visual Gromov hyperbolic space
$X$
the asymptotic dimension of
$X$
is at most the capacity dimension of its boundary at infinity
plus 1,
$\asdim X\le\cdim\di X+1$.
\end{abstract}

\section{Introduction}

The notion of the asymptotic dimension, which is a quasi-isometry
invariant of metric spaces, has been introduced
in \cite{Gr}. The present paper arose as an attempt to fill in
details of a sketch of the proof given in \cite[1.$\text{E}_1'$]{Gr}
that the asymptotic dimension of a negatively pinched
Hadamard manifold
$X$
is bounded above by
$\dim X$, $\asdim X\le\dim X$.
In that way, we came to the notion of the capacity
dimension of a metric space,
$\cdim$,
which should play, as we expect, an important role in
many questions.

Recall that for every Gromov hyperbolic space
$X$
there is a canonical class of metrics on the boundary at infinity
$\di X$
called {\em visual} metrics, see Sect.~\ref{sect:hypspaces}.
Our main result is the following.

\begin{thm}\label{thm:main} Let
$X$
be a visual Gromov hyperbolic space. Then
$$\asdim X\le\cdim\di X+1,$$
for any visual metric on
$\di X$.
\end{thm}

The notion of a visual hyperbolic space (\cite{BoS})
is a rough version of the property that given a base point
$o\in X$,
for every
$x\in X$
there is a geodesic ray in
$X$
emanating from
$o$
and passing through
$x$,
see Sect.~\ref{sect:hypspaces}.

The inequality of Theorem~\ref{thm:main} is sharp.
It is known that
$\asdim\hyp^n=n$
for the real hyperbolic space
$\hyp^n$, $n\ge 2$.
On the other hand, the standard unit sphere metric
is a visual metric on the boundary at infinity
$\di\hyp^n=S^{n-1}$,
and
$\cdim S^{n-1}=n-1$,
see Corollary~\ref{cor:compriem}.

By the definition, the capacity dimension is a bilipschitz
invariant. A remarkable fact
discovered in \cite{LS} is that the close notion of
the Assouad-Nagata dimension is a quasi-symmetry invariant.
It turns out that the capacity dimension is also a
quasi-symmetry invariant, see Sect.~\ref{sect:qsinv},
in particular, the right hand side of the inequality from
Theorem~\ref{thm:main} is independent of the choice of
a visual metric on
$\di X$.
This is also compatible with the fact that every quasi-isometry
of hyperbolic geodesic spaces induces a quasi-symmetry of
their boundaries at infinity.

Now, we briefly describe the structure of the paper. In
Sect.~\ref{sect:prelim} we collect notions and facts from
the dimension theory needed for the paper. Here we also
recall a definition of the asymptotic dimension, see
Sect.~\ref{subsect:asdim}. In Sect.~\ref{sect:cdim} we
give three definitions of the capacity dimension each
of which is useful in different circumstances and prove
their equivalence. Here we also compare the capacity
dimension with the Assouad-Nagata dimension and obtain
monotonicity of the capacity dimension. In Sect.~\ref{sect:qsinv}
we prove that the capacity dimension is a quasi-symmetry
invariant. The proof is based on ideas from \cite{LS}.
The core of the paper is Sect.~\ref{sect:asdim} where
we recall the notion of the hyperbolic cone over a
bounded metric space
$Z$
and prove the relevant estimate for the asymptotic dimension
of the cone via the capacity dimension of the base
$Z$.
In the last Sect.~\ref{sect:hypspaces} we discuss some facts
from the hyperbolic spaces theory and prove Theorem~\ref{thm:main}.

{\em Acknowledgements.} The author thanks Urs Lang and
Viktor Schroeder for inspiring discussions.

\section{Preliminaries}\label{sect:prelim}

Here we collect some (more or less) known notions and facts from
the dimension theory needed in what follows.

Let
$Z$
be a metric space. For
$U$, $U'\sub Z$
we denote by
$\dist(U,U')$
the distance between
$U$
and
$U'$,
$\dist(U,U')=\inf\set{|uu'|}{$u\in U,\ u'\in U'$}$
where
$|uu'|$
is the distance between
$u$, $u'$.
For
$r>0$
we denote by
$B_r(U)$
the open
$r$-neighborhood
of
$U$, $B_r(U)=\set{z\in Z}{$\dist(z,U)<r$}$,
and by
$\ov B_r(U)$
the closed
$r$-neighborhood
of
$U$, $\ov B_r(U)=\set{z\in Z}{$\dist(z,U)\le r$}$.
We extend these notations over all real
$r$
putting
$B_r(U)=U$
for
$r=0$,
and defining
$B_r(U)$
for
$r<0$
as the complement of the closed
$|r|$-neighborhood
of
$Z\sm U$,
$B_r(U)=Z\sm\ov B_{|r|}(Z\sm U)$.

\subsection{Coverings}\label{subsect:cov}

Given a family
$\cU$
of subsets in a metric space
$Z$,
we put
$\mesh(\cU,z)=\sup\set{\diam U}{$z\in U\in\cU$}$
for every
$z\in Z$,
and
$\mesh(\cU)=\sup\set{\diam U}{$U\in\cU$}$.
Clearly,
$\mesh(\cU)=\sup_{z\in Z}\mesh(\cU,z)$.
In the case
$D=\mesh(\cU)<\infty$
we say that
$\cU$
is
$D$-{\em bounded}.

The {\em multiplicity} of
$\cU$, $m(\cU)$,
is the maximal number of members of
$\cU$
with nonempty intersection. For
$r>0$,
the
$r$-{\em multiplicity}
of
$\cU$, $m_r(\cU)$,
is the multiplicity of the family
$\cU_r$
obtained by taking open
$r$-neighborhoods
of the members of
$\cU$.
So
$m_r(\cU)=m(\cU_r)$.
We say that a family
$\cU$
is {\em disjoint} if
$m(\cU)=1$.

A family
$\cU$
is called a {\em covering} of
$Z$
if
$\cup\set{U}{$U\in\cU$}=Z$.
A covering
$\cU$
is said to be {\em colored} if it is the union
of
$m\ge 1$
disjoint families,
$\cU=\cup_{a\in A}\cU^a$, $|A|=m$.
In this case we also say that
$\cU$
is
$m$-colored.
Clearly, the multiplicity of a
$m$-colored
covering is at most
$m$.

Let
$\cU$
be an open covering of a metric space
$Z$.
Given
$z\in Z$
we let
$L'(\cU,z)=\sup\set{\dist(z,Z\sm U)}{$U\in\cU$}$,
$$L(\cU,z)=\min\{L'(\cU,z),\mesh(\cU,z)\}$$
be the Lebesgue number of
$\cU$
at
$z$
(the auxiliary
$L'(\cU,z)$
might be larger than
$\mesh(\cU,z)$
and even infinite as e.g. in the case
$Z=U$
for some member
$U\in\cU$),
$L(\cU)=\inf_{z\in Z}L(\cU,z)$
be the Lebesgue number of
$\cU$.
We have
$L(\cU,z)\le\mesh(\cU,z)$, $L(\cU)\le\mesh(\cU)$
and for every
$z\in Z$
the open ball
$B_r(z)$
of radius
$r=L(\cU)$
centered at
$z$
is contained in some member of the covering
$\cU$.

We shall make use the following

\begin{lem}\label{lem:insidecov} Let
$\cU$
be an open covering of
$Z$
with
$L(\cU)>0$.
Then for every
$s\in(0,L(\cU))$
the family
$\cU_{-s}=B_{-s}(\cU)$
is still an open covering of
$Z$,
and its
$s$-multiplicity
$m_s(\cU_{-s})\le m(\cU)$.
\end{lem}

\begin{proof} For every
$z\in Z$
the ball
$B_r(z)$, $r=L(\cU)$,
is contained in some
$U\in\cU$.
Then
$z\in B_{-s}(U)$
since
$s<r$,
thus
$\cU_{-s}$
is an open covering of
$Z$.
Furthermore, since
$B_s(B_{-s}(U))\sub U$
for every
$U\sub Z$,
we have
$m_s(\cU_{-s})=m(B_s(\cU_{-s}))\le m(\cU)$.
\end{proof}

A covering
$\cU$
is {\em locally finite,} if for every
$z\in Z$
only finitely many its elements intersect some
neighborhood of
$z$.

One defines the {\em nerve} of
$\cU$
as a simplicial polyhedron whose vertex set is
$\cU$,
and a finite collection of vertices spans a simplex iff
the corresponding covering elements have a nonempty intersection.
Thus its (combinatorial) dimension is
$m(\cU)-1$.

\subsection{Uniform polyhedra}\label{subsect:unipol}

Given a index set
$J$,
we let
$R^J$
be the Euclidean space of functions
$J\to\R$
with finite support, i.e.,
$x\in\R^J$
iff only finitely many coordinates
$x_j=x(j)$
are not zero. The distance
$|xx'|$
is well defined by
$$|xx'|^2=\sum_{j\in J}(x_j-x_j')^2.$$
Let
$\De^J\sub\R^J$
be the standard simplex, i.e.,
$x\in\De^J$
iff
$x_j\ge 0$
for all
$j\in J$
and
$\sum_{j\in J}x_j=1$.

A metric in a simplicial polyhedron
$P$
is said to be {\em uniform} if
$P$
is isometric to a subcomplex of
$\De^J\sub\R^J$
for some index set
$J$.
Every simplex
$\si\sub P$
is then isometric to the standard
$k$-simplex
$\De^k\sub\R^{k+1}$, $k=\dim\si$
(so, for a finite
$J$,
$\dim\De^J=|J|-1$).
For every simplicial polyhedron
$P$
there is the canonical embedding
$u:P\to\De^J$,
where
$J$
is the vertex set of
$P$,
which is affine on every simplex. Its image
$P'=u(P)$
is called the {\em uniformization} of
$P$,
and
$u$
is the uniformization map.

For example, the nerve
$\cN=\cN(\cU)$
of a covering
$\cU=\{U_j\}_{j\in J}$
can always be considered as subcomplex of
$\De^J$, $\cN\sub\De^J$,
and therefore as a uniform polyhedron.

\subsection{Barycentric maps}\label{subsect:barycent}

Let
$\cU=\{U_j\}_{j\in J}$
be a locally finite open covering of a metric space
$Z$
by bounded sets,
$\cN=\cN(\cU)\sub\De^J$
its nerve. One defines the {\em barycentric map}
$$p:Z\to\cN$$
associated with
$\cU$
as follows. Given
$j\in J$,
we put
$q_j:Z\to\R$, $q_j(z)=\min\{\diam Z,\,\dist(z,Z\sm U_j)\}$.
Since
$\cU$
is open,
$\sum_{j\in J}q_j(z)>0$
for every
$z\in Z$.
Since
$\cU$
is locally finite and its elements are bounded,
$\sum_{j\in J}q_j(z)<\infty$
for every
$z\in Z$.
Now, the map
$p$
is defined by its coordinate functions
$p_j(z)=q_j(z)/\sum_{j\in J}q_j(z)$, $j\in J$.
Clearly, its image lands at the nerve,
$p(Z)\sub\cN$.
Assume in addition that
$L(\cU)\ge d>0$
and that the multiplicity
$m(\cU)=m+1$
is finite. Then it is easy to check (see for instance \cite{BD},
\cite{BS}) that
$p$
is Lipschitz with Lipschitz constant
$$\Lip(p)\le\frac{(m+2)^2}{d}.$$
Furthermore, for each vertex
$v\in\cN$
the preimage of its open star,
$p^{-1}(\st_v)\sub Z$,
coincides with the member of the covering
$\cU$
corresponding to
$v$.

An (open, locally finite) covering
$\cU'$
is {\em inscribed} in
$\cU$
if every its element is contained in some element of
$\cU$.
In this case there is a simplicial map
$\rho:\cN'\to\cN$
of the nerves which associates to every vertex
$v'\in\cN'$
some vertex
$v\in\cN$
with
$v'\sub v$
(as coverings elements). One easily checks that
this rule is compatible with simplicial structures of
$\cN$, $\cN'$,
and moreover
$\rho\circ p'(z)$
lies in a face of the minimal simplex containing
$p(z)\in\cN$
for every
$z\in Z$.

Note that if
$\mesh(\cU')<L(\cU)$
then
$\cU'$
is inscribed in
$\cU$.

\subsection{Asymptotic dimension}\label{subsect:asdim}

The asymptotic dimension is a quasi-isometry invariant of
a metric space introduced in \cite{Gr}.
There are several equivalent definitions,
see \cite{Gr}, \cite{BD}, and we shall use the following one.
The {\em asymptotic dimension} of a metric space
$X$, $\asdim X$,
is a minimal
$n$
such that for every
$\la>0$
there is a
$\la$-Lipschitz
map
$f:X\to P$
into a uniform simplicial polyhedron
$P$
of dimension
$\le n$
for which the preimages
$f^{-1}(\si)\sub X$
of all simplices
$\si\sub P$
are uniformly bounded. We say that
$f$
is {\em uniformly cobounded} if the last property is satisfied.

\section{Capacity dimension}\label{sect:cdim}

We give three equivalent definitions of the capacity dimension.
Each of them is useful in appropriate circumstances.

Let
$\cU$
be an open covering of a metric space
$Z$.
We define the {\em capacity} of
$\cU$
by
$$\cp(\cU)=\frac{L(\cU)}
         {\mesh(\cU)}\in[0,1],$$
in the case
$\mesh(\cU)=0$
or
$L(\cU)=\mesh(\cU)=\infty$
we put
$\cp(\cU)=1$
by definition.

\subsection{First definition}

For
$\tau>0$, $\de\in(0,1)$
and an integer
$m\ge 0$
we put
$$c_{1,\tau}(Z,m,\de)=\sup_{\cU}\cp(\cU),$$
where the supremum is taken over all open,
$(m+1)$-colored
coverings
$\cU$
of
$Z$
with
$\de\tau\le\mesh(\cU)\le\tau$.

Next, we put
$$c_1(Z,m,\de)=\liminf_{\tau\to 0}c_{1,\tau}(Z,m,\de).$$
The function
$c_1(Z,m,\de)$
is monotone in
$\de$, $c_1(Z,m,\de')\ge c_1(Z,m,\de)$
for
$\de'<\de$.
Hence, there exists a limit
$c_1(Z,m)=\lim_{\de\to 0}c_1(Z,m,\de)$.
Now, we define the {\em capacity dimension} of
$Z$
as
$$\cdim_1(Z)=\inf\set{m}{$c_1(Z,m)>0$}.$$

\subsection{Second definition}

For
$\tau>0$, $\de\in(0,1)$
and an integer
$m\ge 0$
we put
$$c_{2,\tau}(Z,m,\de)=\sup_{\cU}\cp(\cU),$$
where the supremum is taken over all open coverings
$\cU$
of
$Z$
with multiplicity
$\le m+1$
and
$\de\tau\le\mesh(\cU)\le\tau$.
Now, we proceed as above putting
$$c_2(Z,m,\de)=\liminf_{\tau\to 0}c_{2,\tau}(Z,m,\de),$$
$c_2(Z,m)=\lim_{\de\to 0}c_2(Z,m,\de)$
and finally
$$\cdim_2(Z)=\inf\set{m}{$c_2(Z,m)>0$}.$$

\subsection{Third definition}

Let
$f:Z\to P$
be a map into a
$m$-dimensional
uniform polyhedron
$P$.
We define
$\mesh(f)$
as the supremum of
$\diam f^{-1}(\st_v)$
over open stars
$\st_v\sub P$
of vertices
$v\in P$.
Next, we introduce the {\em capacity} of
$f$
as
$$\cp(f)=\left(\Lip(f)\cdot\mesh(f)\right)^{-1},$$
and for
$\tau>0$, $\de\in(0,1)$
and an integer
$m\ge 0$
we define
$c_{3,\tau}(Z,m,\de)=\sup_f\cp(f)$,
where the supremum is taken over all Lipschitz maps
$f:Z\to P$
into
$m$-dimensional
uniform polyhedra
$P$
with
$\de\tau\le\mesh(f)\le\tau$.
Then we define as above
$c_3(Z,m,\de)$, $c_3(Z,m)$
and
$$\cdim_3(Z)=\inf\set{m}{$c_3(Z,m)>0$}.$$

The basic motivation of the capacity dimension is that in some
circumstances we need control over the Lebesgue number
$L(\cU)$
of coverings involved in the definition of a dimension, e.g.,
that the capacity
$\cp(\cU)$
stays separated from zero for appropriately chosen
$\cU$'s.
In general, there is no reason for that. However, if we allow
coverings with larger multiplicity, we can typically gain control
over
$L(\cU)$,
and it may happen that
$\dim Z<\cdim Z$.

Another feature of the definitions is that they involve the auxiliary
variable
$\de$
and the functions
$c_{i,\tau}(Z,m,\de)$, $i=1,2,3$.
This is done for a technical reason to enable ``\v{C}ech
approximations''
$\cU_k$
for which
$\mesh(\cU_{k+1})$
is not extremely small compared with
$\mesh(\cU_k)$
for every
$k$.

\subsection{Equivalence of the definitions}

The proof that three capacity dimensions coincide is standard,
cf.~\cite{BD}, \cite{BS}, \cite{LS}.

\begin{pro}\label{pro:dimcoincide} All three capacity dimensions
coincide,
$$\cdim_1=\cdim_2=\cdim_3.$$
\end{pro}

\begin{proof} The multiplicity of every
$(m+1)$-colored
covering is at most
$m+1$.
Thus
$c_{1,\tau}(Z,m,\de)\le c_{2,\tau}(Z,m,\de)$
for all
$\tau>0$, $\de\in(0,1)$
and integer
$m\ge 0$,
and hence
$\cdim_2(Z)\le\cdim_1(Z)$.

Given an integer
$m\ge 0$,
every open covering
$\cU$
of
$Z$
with multiplicity
$\le m+1$
is locally finite. If in addition its Lebesgue
number is positive,
$L(\cU)>0$,
then the barycentric map
$p:Z\to\cN$,
$\cN=\cN(\cU)$
is the nerve,
$\dim\cN\le m$,
is Lipschitz with
$\Lip(p)\le(m+2)^2/L(\cU)$.
Since
$p^{-1}(\st_v)=U$
for the vertex
$v\in\cN$,
corresponding to
$U\in\cU$,
we have
$\mesh(p)=\mesh(\cU)$.
Thus for capacities we have
$(m+2)^2\cp(p)\ge\cp(\cU)$,
and
$c_{2,\tau}(Z,m,\de)\le(m+2)^2c_{3,\tau}(Z,m,\de)$
for all
$\tau>0$, $\de\in(0,1)$
and integer
$m\ge 0$.
Hence,
$\cdim_3(Z)\le\cdim_2(Z)$.

Finally, we can assume that
$m=\cdim_3(Z)<\infty$.
Then
$c_0=\frac{1}{8}c_3(Z,m)>0$
and
$c_3(Z,m,\de)\ge 4c_0$
for all sufficiently small
$\de>0$, $c_{3,\tau}(Z,m,\de)\ge 2c_0$
for all sufficiently small
$\tau>0$.
Thus there is a Lipschitz map
$f:Z\to P$
into
$m$-dimensional
uniform polyhedron
$P$
with
$\de\tau\le\mesh(f)\le\tau$
and
$\cp(f)\ge c_0$.
Every
$a$-dimensional
simplex
$\si\sub P$
is matched by its barycenter, which is the vertex
$v_\si$
in the first barycentric subdivision
$\ba P$
of
$P$.
We let
$A=\{0,\dots,m\}$
be the color set, and for
$a\in A$
let
$\cP^a$
be the family of the open starts
$\st_{\si}$
of
$\ba P$
of all
$v_\si$
with
$\dim\si=a$,
$\cP^a=\set{\st_{\si}}{$\si\sub P,\ \dim\si=a$}$.
The family
$\cP^a$
is disjoint,
$\st_{\si}\cap\st_{\si'}=\es$
for
$\si\neq\si'$,
and
$\cP=\cup_{a\in A}\cP^a$
is a covering of the polyhedron
$P$.
Thus the open covering
$\cP$
is
$(m+1)$-colored.

Now,
$\cU=p^{-1}(\cP)$
is a
$(m+1)$-colored
open covering of
$Z$, $\cU=\cup_{a\in A}\cU^a$,
where
$\cU^a=p^{-1}(\cP^a)$.
Since the stars
$\st_{\si}$
of
$\ba P$
are contained in appropriate open stars of
$P$,
we have
$\mesh(\cU)\le\mesh(f)\le\tau$.
Since the polyhedron
$P$
is uniform, there is a lower bound
$l_m>0$
for the Lebesgue number of the covering
$\cP$.
Therefore, for the Lebesgue number of
$\cU$
we have
$L(\cU)\ge l_m/\Lip(f)$.
This implies
$\cp(\cU)\ge l_m\cp(f)\ge l_mc_0$
and
$$\mesh(\cU)\ge L(\cU)\ge l_m\cp(f)\mesh(f)
  \ge l_mc_0\de\tau.$$
Putting everything together we obtain
$c_{1,\tau}(Z,m,l_mc_0\de)\ge c_0>0$
for every sufficiently small positive
$\tau$, $\de$.
Thus
$\cdim_1(Z)\le m=\cdim_3(Z)$.
\end{proof}

From now on, we denote by
$\cdim Z$
the common value of the capacity dimensions of
$Z$.
Clearly, the capacity dimension dominates the topological
dimension,
$\dim Z\le\cdim Z$.
The condition for coverings to be open in the first
and second definitions is inessential, and one can
define
$\cdim Z$
using coverings by arbitrary sets.

The following characterization of the capacity dimension
allows to compare it with the Assouad-Nagata dimension,
see \cite{As}, \cite{LS}.

\begin{pro}\label{pro:capass} The capacity dimension of
a metric space
$Z$
is the infimum of all integers
$m$
with the following property: There exists a constant
$c>0$
such that for all sufficiently small
$s>0$,
$Z$
has a
$cs$-bounded
covering with
$s$-multiplicity
at most
$m+1$.
\end{pro}

\begin{proof} We have to prove that
$\cdim Z=\cdim'Z$,
where
$\cdim'Z$
is defined by the statement of the Proposition,
$$\cdim'Z+1=\lim_{c\to\infty}\limsup_{s\to 0}\inf
            \set{m_s(\cU)}{$\mesh(\cU)\le cs$}.$$
Let
$m'=\cdim'Z$.
Then there are positive
$c$, $s_0$
such that for every
$s\in (0,s_0]$
there is a covering
$\cU$
of
$Z$
with
$m_s(\cU)\le m'+1$
and
$\mesh(\cU)\le cs$.
Given a covering
$\cU$
of
$Z$
with
$\mesh(\cU)\le cs$,
note that the covering
$\cU_s=B_s(\cU)$
is open,
$m(\cU_s)=m_s(\cU)$,
and
$$s\le L(\cU_s)\le\mesh(\cU_s)\le\mesh(\cU)+2s\le(c+2)s.$$
Thus for the capacity of
$\cU_s$
we have
$\cp(\cU_s)\ge 1/(c+2)$.
It follows that
$c_{2,\tau}(Z,m',\de)\ge 1/(c+2)$
for
$\tau=(c+2)s$
and
$\de\in (0,1/(c+2))$.
Hence,
$c_2(Z,m')\ge 1/(c+2)>0$
and therefore
$\cdim Z\le m'$.

Conversely, put
$m=\cdim Z$.
Then
$c_0=\frac{1}{8}c_2(Z,m)>0$, $c_2(Z,m,\de)\ge 4c_0$
for all sufficiently small
$\de>0$, $c_{2,\tau}(Z,m,\de)\ge 2c_0$
for all
$\tau$, $0<\tau\le\tau_0$.
This means that there is an open covering
$\cU$
of
$Z$
with
$\de\tau\le\mesh(\cU)\le\tau$
and
$L(\cU)\ge c_0\mesh(\cU)$
having the multiplicity
$\le m+1$.
Then
$s=c_0\de\tau/2<L(\cU)$.
By Lemma~\ref{lem:insidecov} the family
$\cU_{-s}=B_{-s}(\cU)$
is a covering of
$Z$
with
$\mesh(\cU_{-s})\le\mesh(\cU)\le\frac{2}{c_0\de}s$,
and its
$s$-multiplicity
$m_s(\cU_{-s})\le m(\cU)\le m+1$.
Fixing a sufficiently small
$\de$
as above and taking
$\tau\to 0$,
we obtain
$\cdim'Z\le m$.
\end{proof}

From this characterization we immediately obtain

\begin{cor}\label{cor:monotone} The capacity dimension is monotone,
$\cdim Y\le\cdim Z$
for every
$Y\sub Z$.
\qed
\end{cor}

If we omit ``sufficiently small'' from the statement
of Proposition~\ref{pro:capass}, then we come up
with the Assouad-Nagata dimension of
$Z$, $\andim Z$.
Thus we obtain

\begin{cor}\label{cor:capcompareassouad} For every metric space
$Z$
we have
$\cdim Z\le\andim Z$.
\qed
\end{cor}

From \cite{LS} we obtain

\begin{cor}\label{cor:compriem} Every compact Riemannian
manifold
$M$
satisfies
$\cdim M=\dim M$.
\end{cor}

Vice versa, the Assouad-Nagata dimension can be characterized by
the formula
$\andim Z=\inf\set{m\ge 0}{$c_2'(Z,m)>0$}$,
where
$$c_2'(Z,m)=\lim_{\de\to 0}\inf_{\tau>0}c_{2,\tau}(Z,m,\de).$$
Speaking loosely, the Assouad-Nagata dimension takes into
account all sca\-les, while the capacity dimension only all
sufficiently small scales as the topological dimension does.

\section{Quasi-symmetry invariance}\label{sect:qsinv}

The capacity dimension as well as the Assouad-Nagata dimension
is obviously a bilipschitz invariant. The striking fact discovered
in \cite{LS} is that the Assouad-Nagata dimension is a quasi-symmetry
invariant.

A map
$f:X\to Y$
between metric spaces is called {\em quasi-symmetric,}
if it is not constant and there is a homeomorphism
$\eta:[0,\infty)\to[0,\infty)$
such that from
$|xa|\le t|xb|$
it follows that
$|f(x)f(a)|\le\eta(t)|f(x)f(b)|$
for any
$a$, $b$, $x\in X$
and all
$t\ge 0$.
In this case, we say that
$f$
is
$\eta$-{\em quasi-symmetric.}
One easily sees that
$f$
is injective and continuous, and
$f^{-1}:f(X)\to X$
is
$\eta'$-quasi-symmetric
with
$\eta'(t)=1/\eta^{-1}(t^{-1})$
for
$t>0$.
Furthermore, if
$f:X\to Y$
and
$g:Y\to Z$
are
$\eta_f$-
and
$\eta_g$-quasi-symmetric
respectively, then
$g\circ f:X\to Z$
is
$(\eta_g\circ\eta_f)$-quasi-symmetric.
A quasi-symmetric homeomorphism is called a {\em quasi-symmetry}.
For more details on quasi-symmetric maps see \cite{He}.

For example, the transformation
$d\mapsto d^p$
of a metric
$d$,
where
$0<p<1$,
is a quasi-symmetry called a {\em snow-flake transformation},
see e.g. \cite{BoS}. Such a transformation can be far away from
being bilipschitz because nontrivial rectifiable paths
w.r.t. the metric
$d$
are nonrectifiable w.r.t. the metric
$d^p$.

\begin{thm}\label{thm:qsinvariance} The capacity
dimension is a quasi-symmetry invariant of metric
spaces.
\end{thm}

Combining with Corollary~\ref{cor:monotone}, we obtain

\begin{cor}\label{cor:qsmonotone} Assume that there is
a quasi-symmetric
$f:X\to Y$.
Then
$\cdim X\le\cdim Y$.
\qed
\end{cor}

Using Proposition~\ref{pro:capass} one can refer to
\cite[Theorem~1.2]{LS}
for the proof of Theorem~\ref{thm:qsinvariance} because
actually the same argument works in our case. However, we give
a different proof of Theorem~\ref{thm:qsinvariance}, based
on ideas from \cite{LS}, as an attempt to understand
very nice and concise arguments from \cite{LS}.

A key ingredient of the proof of Theorem~\ref{thm:qsinvariance}
is the existence of a sequence of coverings established in
Proposition~\ref{pro:seqcov}.

We say that a family
$\cU$
of sets in a space
$X$
is {\em separated}, if
different members
$U$, $U'\in\cU$
are either disjoint,
$U\cap U'=\es$,
or one of them is contained in the other.
Note that if
$\cU$
is separated, then
$B_{-s}(\cU)$
is separated for every
$s\ge 0$.

Let
$\cU$, $\cU'$
be families of sets in
$X$.
We denote by
$\cU\ast\cU'$
the family obtained by taking for every
$U\in\cU$
the union
$V$
of
$U$
and all members
$U'\in\cU'$
which intersect
$U$, $\cU\ast\cU'=\set{V}{$U\in\cU$}$.
It is straightforward to
check that the following is true.

\begin{lem}\label{lem:undisfam} Let
$\cU$, $\wh\cU$
be separated families in
$X$
such that no member of
$\wh\cU$
intersects disjoint members of
$\cU$.
Then the family
$\cV=\cU\ast\wh\cU\cup\wh\cU$
is separated, moreover, if
$\cU$
is disjoint then
$\cU\ast\wh\cU$
is also disjoint.
\qed
\end{lem}

\begin{pro}\label{pro:seqcov} Suppose that
$X$
is a metric space with finite capacity dimension,
$\cdim X\le n$.
Then there are positive constants
$c_0$, $\de$
such that for every sufficiently small
$r>0$
there exists a sequence of open coverings
$\cU_j$
of
$X$, $j\in\N$,
with the following properties
\begin{itemize}
\item[(i)] For every
$j\in\N$
the covering
$\cU_j$
is
$(n+1)$-colored
by one and the same color set
$A$, $\cU_j=\cup_{a\in A}\cU_j^a$, $|A|=n+1$;
\item[(ii)] for every
$j\in\N$
we have
$$\de r^j\le\mesh(\cU_j)\le r^j\quad
  \text{and}\quad
  L(\cU_j)\ge c_0\mesh(\cU_j);$$
\item[(iii)] for every
$i>j$
the covering
$\cU_i$
is inscribed in the covering
$\cU_j$;
\item[(iv)] for every
$a\in A$
the union
$\cU^a=\cup_{j\in\N}\cU_j^a$
is separated.
\end{itemize}
\end{pro}

\begin{proof} The most important property is (iv). Existence
of a sequence of coverings
$\wh\cU_j$, $j\in\N$,
of
$X$
possessing (i)--(iii) is an easy consequence of the condition
$\cdim X\le n$
and the definition of
$\cdim$.
Namely, as in the proof of Proposition~\ref{pro:capass} we have
$c_0'=\frac{1}{8}c_1(X,n)>0$, $c_1(X,n,\de')\ge 4c_0'$
for all sufficiently small
$\de'>0$.
We fix such a
$\de'$,
and note that
$c_{1,\tau}(X,n,\de')\ge 2c_0'$
for all
$\tau$, $0<\tau\le\tau_0$.
This means that for every
$\tau\in(0,\tau_0]$
there is an open
$(n+1)$-colored
covering
$\cU_\tau$
of
$X$
with
$\de'\tau\le\mesh(\cU_\tau)\le\tau$
and the capacity arbitrarily close to
$c_{1,\tau}(X,n,\de')$,
in particular,
$L(\cU_\tau)\ge c_0'\mesh(\cU_\tau)$.

Take a positive
$r<\min\{c_0'\de'/4,\tau_0\}$
and for every
$j\in\N$
consider the covering
$\wh\cU_j=\cU_{\tau_j}$,
where
$\tau_j=r^j$.
Then the sequence
$\wh\cU_j$, $j\in\N$,
of open coverings of
$X$
satisfies the conditions (i), (ii) (with
$\de=\de'$
and
$c_0=c_0')$.
Since
$L(\wh\cU_j)\ge c_0'\de'r^j>r^{j+1}\ge\mesh(\wh\cU_{j+1})$,
the covering
$\wh\cU_{j+1}$
is inscribed in
$\wh\cU_j$,
thus (iii) is also satisfied.

Now, we modify the sequence
$\{\wh\cU_j\}$
in a way to keep properties (i)--(iii) preserved
and to obtain (iv). This can be done, e.g. as in
\cite[Proposition 4.1]{LS}. We take another track
to show a different possibility.

For
$k\in\N$
we put
$s_k=c_0'\de'r^k/4$.
Fix a color
$a\in A$
and define
$\cV_1^a=\cU_{1,1}^a:=\wh\cU_1^a$.
Assume that for
$k\ge 1$
the family
$\cV_k^a$
is already defined, it is separated and
$\cV_k^a=\cup_{j=1}^k\cU_{j,k}^a$
where each family
$\cU_{j,k}^a$
is disjoint. Now, we define
$$\cV_{k+1}^a:=B_{-s_k}(\cV_k^a)\ast\wh\cU_{k+1}^a\cup\wh\cU_{k+1}^a.$$
Then
$\cV_{k+1}^a=\cup_{j=1}^{k+1}\cU_{j,k+1}^a$,
where
$\cU_{j,k+1}^a=B_{-s_k}(\cU_{j,k}^a)\ast\wh\cU_{k+1}^a$
for
$1\le j\le k$
and
$\cU_{k+1,k+1}^a=\wh\cU_{k+1}^a$.
Since
$\mesh(\wh\cU_{k+1}^a)\le\mesh(\wh\cU_{k+1})\le r^{k+1}<s_k$,
no member of
$\wh\cU_{k+1}^a$
intersects disjoint members of
$B_{-s_k}(\cV_k^a)$.
Then by Lemma~\ref{lem:undisfam} the family
$\cV_{k+1}^a$
is separated and the family
$\cU_{j,k+1}^a$
is disjoint for every
$1\le j\le k+1$.

It follows from the definition of
$\cU_{j,k+1}^a=B_{-s_k}(\cU_{j,k}^a)\ast\wh\cU_{k+1}^a$
and the fact that
$r^{k+1}<s_k$,
that for every
$k\ge 1$, $1\le j\le k$,
every member
$U\in\cU_{j,k+1}^a$
is contained in a unique member
$U'\in\cU_{j,k}^a$.
In this sense, for every
$j\in\N$
the sequence of families
$\cU_{j,k}^a$, $k\ge j$,
is monotone,
$\cU_{j,k}^a\sups\cU_{j,k+1}^a$,
and the intersection
$\cap_{k\ge j}\cU_{j,k}^a$
is defined in the obvious sense.

We put
$\wh s_j=\sum_{k\ge j}s_k=c_0'\de'r^j/4(1-r)$.
Then
$\wh s_j<c_0'\de'r^j\le L(\wh\cU_j)$
for every
$j\ge 1$.
By Lemma~\ref{lem:insidecov}, the family
$\wh\cU_j'=B_{-\wh s_j}(\wh\cU_j)$
is still an open covering of
$X$,
and
$B_{-\wh s_j}(\wh\cU_j^a)\sub\cap_{k\ge j}\cU_{j,k}^a$
for every
$a\in A$, $j\in\N$.
Now, for the interior
$\cU_j^a=\Int\left(\cap_{k\ge j}\cU_{j,k}^a\right)$
the family
$\cU_j=\cup_{a\in A}\cU_j^a$
is an open
$(n+1)$-colored
covering of
$X$
inscribed in
$\wh\cU_j$
for every
$j\in\N$.
Then
$\mesh(\cU_j)\le\mesh(\wh\cU_j)\le r^j$.

Since
$r<1/2$,
we have
$\wh s_j=c_0'\de'r^j/4(1-r)\le c_0'\de'r^j/2$
for every
$j\in\N$.
The covering
$\wh\cU_j'$
is inscribed in
$\cU_j$,
thus
$L(\cU_j)\ge L(\wh\cU_j')\ge c_0'\de'r^j-\wh s_j\ge c_0'\de'r^j/2$.
We put
$\de:=c_0'\de'/2=:c_0$.
Then
$\mesh(\cU_j)\ge L(\cU_j)\ge\de r^j$
and
$L(\cU_j)\ge c_0\mesh(\cU_j)$.
Since
$L(\cU_j)>r^{j+1}\ge\mesh(\cU_{j+1})$,
the covering
$\cU_{j+1}$
is inscribed in
$\cU_j$
for every
$j\in\N$.
Therefore, the sequence of open coverings
$\cU_j$, $j\in\N$,
satisfies (i)--(iii). For every
$a\in A$, $k\ge 1$,
the family
$\cV_k^a$
is separated. It follows that the family
$\cU^a=\cup_{j\in\N}\cU_j^a$
is also separated, hence (iv).
\end{proof}

We define the {\em local capacity} of an open covering
$\cU$
of a metric space
$Z$
by
$$\cploc(\cU)=\inf_{z\in Z}\frac{L(\cU,z)}
         {\mesh(\cU,z)}.$$
Clearly,
$1\ge\cploc(\cU)\ge\cp(\cU)$.
The advantage of the local capacity over the capacity
is that its positivity is preserved under quasi-symmetries
quantitatively, see Lemma~\ref{lem:capacity}.
This implies that a dimension defined exactly as the
capacity dimension replacing the capacity of coverings by
the local capacity is a quasi-symmetry invariant. However,
that invariant is not as good as the capacity dimension for
applications. We use the local capacity of coverings only
as an auxiliary tool to prove quasi-symmetry invariance
of the capacity dimension.

\begin{lem}\label{lem:capacity} Let
$\cU$
be an open covering of a metric space
$Z$, $f:Z\to Z'$
be an
$\eta$-quasi-symmetry,
$\cU'=f(\cU)$
be the image of
$\cU$.
Then for the local capacities of
$\cU$
and
$\cU'$
we have
$$\frac{1}{\cploc(\cU')}\le 16\eta\left(\frac{2}{\cploc(\cU)}\right).$$
\end{lem}

\begin{proof} We can assume that
$\cploc(\cU)>0$.
We fix
$z\in Z$
and consider
$U\in\cU$
for which
$z\in U$
and
$\dist(z,Z\sm U)\ge L(\cU,z)/2$.
For
$z'=f(z)$
and
$U'=f(U)$
there is
$a'\in Z'\sm U'$
with
$|z'a'|\le 2\dist(z',Z'\sm U')$.
Then
$|z'a'|\le 2L(\cU',z')$,
and for
$a=f^{-1}(a')$
we have
$|za|\ge\dist(z,Z\sm U)\ge L(\cU,z)/2$.

Similarly, consider
$V'\in\cU'$
with
$z'\in V'$
and
$\diam V'\ge\mesh(\cU',z')/2$.
There is
$b'\in V'$
with
$|z'b'|\ge\diam V'/4$.
Then
$|z'b'|\ge\mesh(\cU',z')/8$,
and for
$b:=f^{-1}(b')$
we have
$|zb|\le\mesh(\cU,z)$.
Therefore, we have
$$\cploc(\cU)\le\frac{L(\cU,z)}{\mesh(\cU,z)}
  \le\frac{2|za|}{|zb|}\quad\text{and}\quad
  |zb|\le t|za|$$
with
$t=2/\cploc(\cU)$.
It follows
$|z'b'|\le\eta(t)|z'a'|$
and
$$\frac{L(\cU',z')}{\mesh(\cU',z')}\ge\frac{|z'a'|}{16|z'b'|}
  \ge\left(16\eta(t)\right)^{-1}$$
for every
$z'\in Z'$.
Then
$\cploc(\cU')\ge
 \left(16\eta\left(\frac{2}{\cploc(\cU)}\right)\right)^{-1}$.
\end{proof}

A covering
$\cU$
of
$Z$
is said to be
$c$-{\em balanced}, $c>0$,
if
$\inf\set{\diam(U)}{$U\in\cU$}\ge c\cdot\mesh(\cU)$.
The notion of a balanced covering combined with the
local capacity allows to estimate from below the
capacity of a covering as follows.

\begin{lem}\label{lem:balocap} If an open covering
$\cU$
of a metric space
$Z$
is
$c_1$-balanced
and its local capacity satisfies
$\cploc(\cU)\ge c_0$,
then
$\cp(\cU)\ge c_0\cdot c_1$.
\end{lem}

\begin{proof} Since
$\cU$
is
$c_1$-balanced,
we have
$\mesh(\cU,z)\ge c_1\cdot\mesh(\cU)$
for every
$z\in Z$.
Since
$\cploc(\cU)\ge c_0$,
we have
$L(\cU,z)\ge c_0\cdot\mesh(\cU,z)$
for every
$z\in Z$.
Therefore
$L(\cU)\ge c_0c_1\cdot\mesh(\cU)$.
\end{proof}

Let
$f:X\to Y$
be a quasi-symmetry. To prove Theorem~\ref{thm:qsinvariance}
it suffices to show that
$\cdim Y\le\cdim X$.
The idea is to construct out of a sequence
$\cU_j$, $j\in\N$,
of coverings of
$X$
as in Proposition~\ref{pro:seqcov} a covering
$\cV$
of
$X$
with local capacity uniformly separated from 0,
see Lemma~\ref{lem:minloc}, such that its image
$f(\cV)$
has an arbitrarily small mesh and is balanced,
see Lemma~\ref{lem:minbalance}. Then by
Lemma~\ref{lem:balocap} combined with Lemma~\ref{lem:capacity}
the capacity of the covering
$f(\cV)$
of
$Y$
is positive quantitatively, which implies
$\cdim Y\le\cdim X$.

Fix a sufficiently small
$r>0$
and consider the sequence of open coverings
$\cU_j$
of
$X$
as in Proposition~\ref{pro:seqcov}. We can assume additionally that
$\diam U\ge L(\cU_j)$
for every
$U\in\cU_j$,
since if
$\diam U<L(\cU_j)$
then
$U$
is contained in another member
$U'\in\cU_j$
and thus it can be deleted from
$\cU_j$
without destroying any property from (i)--(iv).

Following \cite{LS} we put
$\cU=\cup_{j\in\N}\cU_j$
and for
$s>0$
consider the family
$\cU(s)=\set{U\in\cU}{$\diam f(U)\le s$}$.

\begin{lem}\label{lem:scover} For every
$s>0$
the family
$\cU(s)$
is a covering of
$X$.
\end{lem}

\begin{proof} We fix
$x\in X$,
consider
$x'\in X$
different from
$x$
and put
$y=f(x)$, $y'=f(x')$.
For every
$j\in\N$
there is
$U_j\in\cU_j$
containing
$x$.
Take
$y_j\in f(U_j)$
with
$\diam f(U_j)\le 4|yy_j|$
and consider
$x_j=f^{-1}(y_j)$.
Then
$|xx_j|\le t_j|xx'|$
with
$t_j\to 0$
as
$j\to\infty$,
since
$\diam U_j\le r^j\to 0$.
Therefore,
$\diam f(U_j)\le 4|yy_j|\le 4\eta(t_j)|yy'|\le s$
for sufficiently large
$j$.
Hence,
$U_j\in\cU(s)$.
\end{proof}

A family
$\cV\sub\cU(s)$
is {\em minimal} if every
$U\in\cU(s)$
is contained in some
$V\in\cV$
and neither of different
$V$, $V'\in\cV$
is contained in the other.

\begin{lem}\label{lem:mincover} For every
$s>0$
there is a minimal family
$\cV\sub\cU(s)$.
Every minimal family
$\cV\sub\cU(s)$
is a
$(n+1)$-colored
covering of
$X$.
\end{lem}

\begin{proof} Given
$s>0$
we construct a family
$\cV\sub\cU(s)$
deleting every
$U\in\cU(s)$
which is contained in some other
$U'\in\cU(s)$.
Now,
$\cV$
is what remains. One needs only to check that for
every
$U\in\cU(s)$
there is a maximal
$U'\in\cU(s)$
with
$U\sub U'$.
It follows from Proposition~\ref{pro:seqcov}(i)
that for every
$j\in\N$
there are only finitely many
$U'\in\cU_j$
containing
$U$
(since all of them must have different colors).
Since
$\mesh(\cU_j)\to 0$
as
$j\to\infty$,
there are only finitely many
$U'\in\cU(s)$
containing
$U$
and hence
there is a maximal
$U'\in\cU(s)$
among them.

Let
$\cV\sub\cU(s)$
be a minimal family. By Lemma~\ref{lem:scover}, the family
$\cU(s)$
is a covering of
$X$,
and it follows from the definition of a minimal family that
$\cV$
is also a covering of
$X$.
It follows from Proposition~\ref{pro:seqcov}(iv) that
different
$V$, $V'\in\cV$
having one and the same color are disjoint. Thus
$\cV$
is
$(n+1)$-colored.
\end{proof}

\begin{lem}\label{lem:minloc} There is a constant
$\nu>0$
depending only on
$c_0$, $\de$, $r$
and
$\eta$
such that for every
$s>0$
every minimal covering
$\cV\sub\cU(s)$
has the local capacity
$\cploc(\cV)\ge\nu$.
\end{lem}

\begin{proof} Let
$\cV\sub\cU(s)$
be a minimal family. Given
$x\in X$
we put
$j=j(x)=\min\set{i\in\N}{$x\in V\in\cV\cap\cU_i$}$.
Then
$\mesh(\cV,x)\le r^j$.
We fix
$V\in\cV\cap\cU_j$
containing
$x$, $v\in V$
with
$4|xv|\ge\diam V$
and note that
$\diam V\ge L(\cU_j)\ge c_0\de r^j$
by our assumptions.

Furthermore, we fix
$\mu>0$
with
$4\eta(4\mu/c_0\de)\le 1$.
Now we check that for
$i\in\N$
with
$r^{i-j}\le\mu$
every
$U\in\cU_i$
containing
$x$
is a member of
$\cU(s)$.
There is
$u\in U$
with
$4|f(x)f(u)|\ge\diam f(U)$.
We have
$|xu|\le t|xv|$
for some
$t\le 4\diam U/\diam V\le 4r^{i-j}/c_0\de\le 4\mu/c_0\de$.
Then
$\diam f(U)\le 4|f(x)f(u)|\le 4\eta(4\mu/c_0\de)|f(x)f(v)|
  \le\diam V\le s$,
thus
$U\in\cU(s)$.

Therefore,
$L(\cV,x)\ge L(\cU_i)\ge c_0\de r^i$.
Assuming that
$i$
is taken minimal with
$r^{i-j}\le\mu$,
we obtain
$L(\cV,x)\ge c_0\de\mu r^{j+1}$.
Thus
$\frac{L(\cV,x)}{\mesh(\cV,x)}\ge\nu=c_0\de\mu r$
for every
$x\in X$
and
$\cploc(\cV)\ge\nu$.
\end{proof}

\begin{lem}\label{lem:minbalance} Given a minimal
family
$\cV\sub\cU(s)$,
the
$(n+1)$-colored
covering
$\cW=f(\cV)$
of
$Y$
satisfies
$\diam W\ge s/4\eta(t)$
for every
$W\in\cW$,
where
$t=4/c_0\de r$.
In particular,
$\mesh(\cW)\ge s/4\eta(t)$
and
$\cW$
is
$c$-balanced
with
$c\ge1/4\eta(t)$.
\end{lem}

\begin{proof} Note that
$\mesh(\cW)\le s$
by the definition of
$\cU(s)$.
Take any
$W\in\cW$
and consider
$V=f^{-1}(W)$.
We can assume that
$V\in\cU_j$
for some
$j\in\N$.
Then
$\diam V\ge L(\cU_j)\ge c_0\de r^j$
by our assumption on the sequence
$\{\cU_j\}$.

For any
$U\in\cU$
with
$V\sub U$
we have
$\diam f(U)>s$,
since the family
$\cV$
is minimal. The covering
$\cU_j$
is inscribed in
$\cU_{j-1}$,
thus there is
$U\in\cU_{j-1}$
containing
$V$,
in particular,
$\diam f(U)>s$.

Take
$y\in W\sub f(U)$.
There is
$y'\in f(U)$
with
$|yy'|\ge\diam f(U)/4>s/4$.
For
$x=f^{-1}(y)$, $x'=f^{-1}(y')$
we have
$|xx'|\le\diam U\le\mesh(\cU_{j-1})\le r^{j-1}$.
There is
$v\in V$
with
$|xv|\ge\diam V/4\ge c_0\de r^j/4$.
Thus
$|xx'|\le r^{j-1}\le t|xv|$
for
$t=4/c_0\de r$.
For
$w=f(v)\in W$
we obtain
$|yy'|\le\eta(t)|yw|\le\eta(t)\diam W$.
Hence,
$\diam W\ge s/4\eta(t)$.
\end{proof}

\begin{proof}[Proof of Theorem~\ref{thm:qsinvariance}]
Let
$f:X\to Y$
be an
$\eta$-quasi-symmetry.
We show that
$\cdim Y\le n$
for every
$n\ge\cdim X$.
Fix a sufficiently small
$r>0$
and consider a sequence
$\cU_j$, $j\in\N$,
of coverings of
$X$
as in Proposition~\ref{pro:seqcov}
with positive
$c_0$, $\de$.
Then by Lemmas~\ref{lem:mincover} and
\ref{lem:minbalance} for every
$s>0$
we have an open
$(n+1)$-colored
covering
$\cW$
of
$Y$
with
$s/4\eta(t)\le\mesh(\cW)\le s$,
which is
$c$-balanced,
$c\ge 1/4\eta(t)$,
where
$t=4/c_0\de r$.
Moreover, by Lemmas~\ref{lem:minloc} and
\ref{lem:capacity} its local capacity
$\cploc(\cW)\ge d$,
where the constant
$d>0$
depends only on
$\eta$, $c_0$, $\de$, $r$.
Then by Lemma~\ref{lem:balocap} we have
$\cp(\cW)\ge c\cdot d$
independently of
$s$.

This shows that
$c_{1,s}(Y,n,\de')\ge c\cdot d$
for every
$s>0$,
where
$\de'=1/4\eta(t)$.
Hence,
$c_1(Y,n)\ge c\cdot d>0$
and
$\cdim Y\le n$.
\end{proof}

\section{Asymptotic dimension of a hyperbolic cone}\label{sect:asdim}

Let
$Z$
be a bounded metric space. Assuming that
$\diam Z>0$
we put
$\mu=\pi/\diam Z$
and note that
$\mu|zz'|\in[0,\pi]$
for every
$z$, $z'\in Z$.
Recall that the hyperbolic cone
$\cone(Z)$
over
$Z$
is the space
$Z\times[0,\infty)/Z\times\{0\}$
with metric defined as follows. Given
$x=(z,t)$, $x'=(z',t')\in\cone(Z)$
we consider a triangle
$\ov o\,\ov x\,\ov x'\sub\hyp^2$
with
$|\ov o\,\ov x|=t$, $|\ov o\,\ov x'|=t'$
and the angle
$\angle_{\ov o}(\ov x,\ov x')=\mu|zz'|$.
Now, we put
$|xx'|:=|\ov x\,\ov x'|$.
In the degenerate case
$Z=\{\pt\}$
we define
$\cone(Z)=\{\pt\}\times[0,\infty)$
as the metric product. The point
$o=Z\times\{0\}\in\cone(Z)$
is called the {\em vertex} of
$\cone(Z)$.

\begin{thm}\label{thm:conedim} For every bounded metric space
$Z$
we have
$$\asdim\cone(Z)\le\cdim Z+1.$$
\end{thm}

The proof occupies sect.~\ref{subsect:esthyp}--\,\ref{subsect:proof}.
In sect.~\ref{subsect:homotopy}--\,\ref{subsect:proof} our
arguments are close to those from \cite[\S2]{BD}.

\subsection{Some estimates from hyperbolic geometry}\label{subsect:esthyp}

We denote by
$Z_t$
the metric sphere of radius
$t>0$
around
$o$
in
$\cone(Z)$.
There are natural polar coordinates
$x=(z,t)$, $z\in Z$, $t\ge 0$,
in
$\cone(Z)$.
Then
$Z_t=\set{(z,t)}{$z\in Z$}$
is the copy of
$Z$
at the level
$t$.
For
$t>0$
we denote by
$\pi_t:Z_t\to Z$
the canonical homeomorphism,
$\pi_t(z,t)=z$.

Let
$\cU$
be an open covering of
$Z$
with multiplicity
$m+1$
and positive Lebesgue number
$L(\cU)$.
Let
$\cN=\cN(\cU)$
be the nerve of
$\cU$,
$p:Z\to\cN$
the barycentric map. Then
$\Lip(p)\le\frac{(m+2)^2}{L(\cU)}$,
see sect.~\ref{subsect:barycent}.
For every
$t>0$
we have the induced covering
$\cU_t=\pi_t^{-1}(\cU)$
of
$Z_t$
whose nerve is canonically isomorphic to
$\cN$,
and the corresponding barycentric map
$p_t:\cU_t\to\cN$.

Given
$\la>0$
we want to find
$t>0$
and conditions for
$\cU$
such that
$\Lip(p_t)\le\la$
and still to get
$p_t$
uniformly cobounded w.r.t. the metric induced from
$\cone(Z)$.
To this end, we first recall the hyperbolic cosine law. For
$t>0$, $\al\in[0,\pi]$
we define
$a=a(t,\al)$
by
$$\cosh a=\cosh^2(t)-\sinh^2(t)\cos\al,$$
i.e.,
$a$
is the length of the base opposite to the vertex
$o$
with angle
$\al$
of a isosceles triangle in
$\hyp^2$
with sides
$t$.
Then for
$\al$
sufficiently small we have
$$\cosh a=1+\frac{1}{2}\sinh^2(t)\al^2+\sinh^2(t)\cdot o(\al^3).$$
Assume that small
$\la$, $\si>0$
are fixed so that
$d:=\frac{(m+2)^2}{\la}-\ln\frac{1}{\si}=\frac{(m+2)^2}{2\la}$,
and
$$\si\tau\le \mu L(\cU)\le\mu\mesh(\cU)\le\tau$$
for sufficiently small
$\tau$.
We put
$$t_\tau=\ln\frac{2}{\tau}+\frac{2(m+2)^2}{\la}.$$
Then
$t_\tau-2d=\ln\frac{2}{\si^2\tau}$,
and for
$t_\tau-2d\le t\le t_\tau$
we have
$$\sinh^2(t)(\si\tau)^2\simeq\frac{1}{4}e^{2t}(\si\tau)^2
  \ge\frac{1}{\si^2}=\exp\left(\frac{(m+2)^2}{\la}\right)\gg 1,$$
while
$$\sinh^2(t)\cdot o(\tau^3)
  \le o(\tau)\cdot\exp\left(\frac{4(m+2)^2}{\la}\right)
  \ll 1.$$
Noting that
$L(\cU_t)=a(t,\mu L(\cU))$
we obtain
$$\cosh(L(\cU_t))\ge 1+\frac{1}{2}\sinh^2(t)(\si\tau)^2
  \ge 1+\frac{1}{2}
  \exp\left(\frac{(m+2)^2}{\la}\right)$$
up to a negligible error. Hence
$L(\cU_t)>\frac{(m+2)^2}{\la}$
and
$\Lip(p_t)\le\frac{(m+2)^2}{L(\cU_t)}<\la$.

Similarly,
$\mesh(\cU_t)=a(t,\mu\mesh(\cU))$,
and for
$t_\tau-2d\le t\le t_\tau$,
we obtain as above
$$\cosh(\mesh(\cU_t))\le 1+\frac{1}{2}\sinh^2(t)\tau^2
  \simeq 1+\frac{1}{8}
   e^{2t}\tau^2\le 1+\frac{1}{2}
  \exp\left(\frac{4(m+2)^2}{\la}\right),$$
which gives an upper bound for
$\mesh(\cU_t)$
depending only on
$\la$.

\subsection{\v{C}ech approximation}\label{subsect:cech}

Here we construct a sequence of coverings
$\{\cU_k\}$
and associated barycentric maps
which will be used in the proof of Theorem~\ref{thm:conedim}.
We can assume that the capacity dimension of
$Z$
is finite,
$m=\cdim Z<\infty$.
Then
$c_0=\frac{1}{8}c_2(Z,m)>0$,
and
$c_2(Z,m,\de)\ge 4c_0$
for all sufficiently small
$\de>0$,
see sect.~\ref{sect:cdim}.
Given
$\la>0$
we take
$\de>0$
so that
$$d:=\frac{(m+2)^2}{\la}+\ln(c_0\de)=\frac{(m+2)^2}{2\la},$$
assuming that
$\la$
is sufficiently small to satisfy
$c_2(Z,m,\de)\ge 4c_0$.
Then
$c_{2,\tau}(Z,m,\de)\ge 2c_0$
for all
$\tau$, $0<\tau\le\tau_0$.
Consider the sequence
$\tau_k$, $k\ge 0$,
recursively defined by
$\tau_{k+1}=e^{-2d}\tau_k$.
Recall that
$\mu=\pi/\diam Z$.
Then it follows from the second definition of
$\cdim Z$
that for every
$k\ge 0$
there is an open covering
$\wh\cU_k$
of
$Z$
such that

\begin{itemize}

\item[(i)] $m(\wh\cU_k)\le m+1$;
\item[(ii)] $\de\tau_k\le\mu\mesh(\wh\cU_k)\le\tau_k$
and $L(\wh\cU_k)\ge c_0\mesh(\wh\cU_k)$.

\end{itemize}
Since
$\mu L(\wh\cU_k)\ge c_0\de\tau_k>\tau_{k+1}\ge\mu\mesh(\wh\cU_{k+1})$,
we additionally have

\begin{itemize}

\item[(iii)] $\wh\cU_{k+1}$
is inscribed in
$\wh\cU_k$
for every
$k$,
\end{itemize}
cf. Proposition~\ref{pro:seqcov}. We put
$$t_k=\ln\frac{2}{\tau_k}+\frac{2(m+2)^2}{\la},$$
$Z_k=Z_{t_k}\sub\cone(Z)$
and for every
$t>0$
and every integer
$k\ge 0$
consider the covering
$\cU_{t,k}=\pi_{t}^{-1}(\wh\cU_k)$
of
$Z_t$.
Note that its nerve is independent of
$t>0$
and can be identified with
$\cN_k=\cN(\wh\cU_k)$.
Then
$t_k-t_{k-1}=2d$,
and using the estimates from sect.~\ref{subsect:esthyp}
with
$\tau=\tau_k$
and
$\si=c_0\de$
we obtain for the barycentric map
$p_{t,k}:Z_t\to\cN_k$
associated with the covering
$\cU_{t,k}$
that
$\Lip(p_{t,k})<\la$
for all
$t_{k-1}\le t\le t_k$
and all
$k\ge 1$.
We put
$\cU_k=\cU_{t_k,k}$
and
$p_k=p_{t_k,k}:Z_k\to\cN_k$.
Then
$\Lip(p_k)<\la$.
Furthermore,
$\mesh(\cU_k)$
is bounded above by a constant depending only on
$\la$, $\mesh(\cU_k)\le\const(\la)$,
for all
$k$.
Hence, preimages of all simplices from
$\cN_k$
under
$p_k$
have uniformly bounded diameter
$\le\const(\la)$
independently of
$k$.

\subsection{Homotopy between
$p_k$
and
$\rho_k\circ p_{k+1}$}\label{subsect:homotopy}

Due to (iii), for every
$k$
there is a simplicial map
$\rho_k:\cN_{k+1}\to\cN_k$
such that
$\rho_k\circ p_{k+1}(z,t_{k+1})$
lies in a face of the minimal simplex containing
$p_k(z,t_k)\in\cN_k$
for every
$z\in Z$.

\begin{lem}\label{lem:lipsimplicial} For every
$k$,
the map
$\rho_k:\cN_{k+1}\to\cN_k$
is
$c_1$-Lipschitz
with
$c_1=c_1(m)$
depending only on
$m$.
\end{lem}

\begin{proof} Recall that the nerve
$\cN_k$
is a uniform polyhedron, and that
$\rho_k$
is affine on every simplex sending it
either to an isometric copy in
$\cN_k$
or to a face of it. In either case,
$\rho_k$
is
$c_1$-Lipschitz
on every simplex.
If
$x$, $x'\in\cN_{k+1}$
are from disjoint simplices, then
$|x-x'|\ge\frac{\sqrt 2}{c_1}$
for some
$c_1=c_1(m)>0$,
and
$|\rho_k(x)-\rho_k(x')|\le\sqrt 2$.
The remaining case, when
$x$, $x'$
are sitting in different simplices in
$\cN_{k+1}$
having a common face, we leave to the reader
as an exercise.
\end{proof}

Consider the annulus
$A_k\sub\cone(Z)$
between
$Z_k$
and
$Z_{t_k+d}$,
$A_k=Z\times[t_k,t_k+d]$
in the polar coordinates. We put
$s=s_k(t)=\frac{1}{d}(t-t_k)$
for
$t_k\le t\le t_k+d$
and define the homotopy
$h_k:A_k\to\cN_k\times[0,1]$
between
$p_k$
and
$\rho_k\circ p_{k+1}$
by
$$h_k(z,t)=\left((1-s)p_k(z,t_k)
  +s\rho_k\circ p_{k+1}(z,t_{k+1}),s\right).$$
This is well defined because the points
$p_k(z,t_k)$
and
$\rho_k\circ p_{k+1}(z,t_{k+1})$, $z\in Z$,
can be joined by the segment in the appropriate simplex.

\begin{lem}\label{lem:liphomotopy} The map
$h_k$
is
$c_2\la$-Lipschitz
with respect to the product metric on
$\cN_k\times[0,1]$,
for some constant
$c_2=c_2(m)>0$
depending only on
$m$.
\end{lem}

\begin{proof} By convexity of the distance function
in Euclidean space, the distance between
$h_k(z,t)$, $h_k(z',t)\in\cN_k\times\{s\}$
is bounded above by the maximum of distances between the
end points of the vertical segments
$z\times[0,1]$, $z'\times[0,1]$
containing them, thus
$$|h_k(x_t)-h_k(x_t')|\le\max\{|p_k(x_k)-p_k(x_k')|,
  |\rho_k\circ p_{k+1}(x_{k+1})-\rho_k\circ p_{k+1}(x_{k+1}')|\},$$
where
$t_k\le t\le t_k+d$, $x_t=(z,t)$, $x_t'=(z',t)$,
$x_k=x_{t_k}$, $x_k'=x_{t_k}'$.
On the other hand,
$$|p_k(x_k)-p_k(x_k')|\le\Lip(p_k)|x_kx_k'|
 \le\la|x_tx_t'|,$$
and by Lemma~\ref{lem:lipsimplicial},

\begin{eqnarray*}
  |\rho_k\circ p_{k+1}(x_{k+1 })-\rho_k\circ p_{k+1}(x_{k+1}')|
  &=&|\rho_k\circ p_{t,k+1}(x_t)-\rho_k\circ p_{t,k+1}(x_t')|\\
  &\le& c_1\Lip(p_{t,k+1})|x_tx_t'|\le c_1\la|x_tx_t'|
\end{eqnarray*}
for
$t_k\le t\le t_k+d$.
Furthermore, since every edge of any standard simplex has
length
$\sqrt 2$,
we have
$$|h_k(z',t)-h_k(z',t')|\le\sqrt 3|s-s'|=
  \frac{\sqrt 3}{d}|(z',t)(z',t')|,$$
where
$s'=s_k(t')$.
Taking into account that
$$|(z,t)(z',t')|\ge\max\{|x_tx_t'|,|(z',t)(z',t')|\}$$
(assuming
$t'\ge t$),
and that
$d=\frac{(m+2)^2}{2\la}$,
we obtain from all of these that
$\Lip(h_k)\le c_2\la$.
\end{proof}

\subsection{Simplicial mapping cylinder of
$\rho_k$}\label{subsect:mapcyl}

Consider the annulus
$B_k\sub\cone(Z)$
between
$Z_{t_k+d}$
and
$Z_{k+1}$, $B_k=Z\times[t_k+d,t_{k+1}]$
in the polar coordinates, and define
$g_k:B_k\to\cN_{k+1}\times[0,1]$
as follows
$g_k(z,t_k+d)=(p_{k+1}(z,t_{k+1}),0)$,
$g_k(z,t_{k+1})=(p_{k+1}(z,t_{k+1}),1)$
and
$g_k$
is affine on every segment
$z\times[t_k+d,t_{k+1}]\sub B_k$, $z\in Z$.

Since
$\frac{1}{d}<\la$, it immediately follows from the
estimates of sect.~\ref{subsect:cech} that
$g_k$
is
$\la$-Lipschitz
(with respect to the product metric on
$\cN_{k+1}\times[0,1]$).

Next, recall the notion of the simplicial mapping cylinder
for a simplicial map
$\rho:K\to L$
of simplicial complexes (see, e.g., \cite{Sp}). Assuming that the
vertices of
$K$
are linearly ordered, we define the mapping cylinder
$C_\rho$
of
$\rho$
as a simplicial complex whose vertex set is the union of
those of
$K$
and
$L$,
and simplices are the simplices of
$K$
and
$L$
and all subsets of the sets
$\{v_0,\dots,v_k,\rho(v_k),\dots,\rho(v_p)\}$,
where
$\{v_0<\dots<v_p\}$
is a simplex in
$K$.

Now, assuming that a linear order on
$\cN_{k+1}$
is fixed, one triangulates
$\cN_{k+1}\times[0,1]$
as the mapping cylinder of the identity map, and defines
the canonical simplicial map
$\phi_k:\cN_{k+1}\times[0,1]\to C_k=C_{\rho_k}$,
which sends
$\cN_{k+1}\times\{0\}$
onto the subcomplex
$\rho_k(\cN_{k+1})\sub C_k$
by
$\rho_k$,
and
$\cN_{k+1}\times\{1\}$
onto
$\cN_{k+1}\sub C_k$
identically.

By \cite[Proposition~3]{BD},
$\phi_k$
is
$c(m)$-Lipschitz
for some constant
$c(m)>0$
depending only on
$m\ge\dim\cN_{k+1}$,
where the cylinder
$C_k$
is given the uniform metric (there is a minor
inaccuracy in argument there claiming that
a simplicial map between uniform complexes is
always 1-Lipschitz, which is not true as easily seen
for
$\De^m\to\De^1$
with
$m\ge 2$;
it is only
$c(m)$-Lipschitz,
see Lemma~\ref{lem:lipsimplicial}). In conclusion, the composition
$\phi_k\circ g_k:B_k\to C_k$
is
$c_3\la$-Lipschitz
with
$c_3=c_3(m)$.
Note that
$h_k$
and
$\phi_k\circ g_k$
coincide on
$Z_{t_k+d}=A_k\cap B_k\sub\cone(Z)$
if one identifies
$\rho_k(\cN_{k+1})\sub C_k$
with subcomplex in
$\cN_k\times\{1\}$.

\subsection{Proof of Theorem~\ref{thm:conedim}}\label{subsect:proof}

We have to find for every sufficiently small
$\la>0$
a uniform polyhedron
$P$
with
$\dim P\le\cdim Z+1$
and a uniformly cobounded,
$\la$-Lipschitz
map
$f:\cone(Z)\to P$.

Given
$\la>0$
we take
$\de>0$
as in sect.~\ref{subsect:cech}. Then we generate
sequences
$\{\tau_k\}$
of positive reals,
$\{\cU_k\}$
of open coverings of
$Z$, $k\ge 0$,
and all the machinery around them from
sect.~\ref{subsect:cech}--\,\ref{subsect:mapcyl}.

Having that we define
$P$
as the uniformization of the union
$P'=P_{-1}\cup_{k\ge 0}P_k$,
where
$P_k$
is constructed out of the uniformization of
$\cN_k\times[0,1]$
(triangulated by fixing a linear order on
$\cN_k$)
and the simplicial mapping cylinder
$C_k$
by attaching them along the common subcomplex
$\rho_k(\cN_{k+1})\sub(\cN_k\times\{1\})\cap C_k$.
Furthermore,
$P_{k+1}$
is attached to
$P_k$
along the common subcomplex
$\cN_{k+1}$
for every
$k\ge 0$.
The polyhedron
$P_{-1}$
is the cone over
$\cN_0$
attached to
$P_0$
along the base. Then
$\dim P\le m+1$
for
$m=\cdim Z$.

The map
$f:\cone(Z)\to P$
is obtained by composing the map
$f':\cone(Z)\to P'$
with the uniformization of
$P'$,
where
$f'$
coincides with
$h_k$
on
$A_k$
and with
$\phi_k\circ g_k$
on
$B_k$
for every
$k\ge 0$.
Finally,
$f'$
is affine on every segment
$z\times[0,t_0]$, $z\in Z$,
sending
$o=Z\times\{0\}$
to the vertex of
$P_{-1}$.
It follows from Lemma~\ref{lem:liphomotopy} and
sect.~\ref{subsect:mapcyl} that
$f$
is
$c\la$-Lipschitz
for some
$c=c(m)>0$
on every
$A_k$, $B_k$, $k\ge 0$,
and on
$Z\times[0,t_0]\sub\cone(Z)$.

Since
$\diam P\le\sqrt 2$
and
$\frac{1}{d}<\la$,
the
$c\la$-Lipschitz
condition is certainly satisfied for points
$(z,t)$, $(z',t')\in\cone(Z)$
separated by some annulus
$A_k$
or
$B_k$.
Thus we assume that
$(z,t)$, $(z',t')$
are sitting in adjacent annuli. Unfortunately,
we cannot directly apply the argument from \cite[Proposition~4]{BD}
which would be well adapted to our situation if
$\cone(Z)$
is geodesic. In general, this is not the case, and we
slightly modify it as follows.

Assume that
$t'>t$.
We take
$t''\in(t,t')$
for which
$(z',t'')$
is common for the annuli, and note that
$$|(z,t)(z',t')|\ge\max\{|(z,t)(z',t)|,|t-t'|=|t-t''|+|t''-t'|\}$$
by geometry of
$\cone(Z)$.
Now, the required Lipschitz condition for the pair
$(z,t)$
and
$(z',t')$
follows in the obvious way from those for three pairs
$(z,t)$
and
$(z',t)$,
$(z',t)$
and
$(z',t'')$,
$(z',t'')$
and
$(z',t')$,
each of which belong to some annulus.

It remains to check that
$f$
is uniformly cobounded. For every simplex
$\si\sub P_k$, $k\ge 0$,
the preimage
$f'^{-1}(\si)\sub A_k\cup B_k$
is contained in
$Z_\si\times[t_k,t_{k+1}]\sub\cone(Z)$,
where
$\diam(Z_\si\times\{t_k\})\le\const(\la)$
by estimates of sect.~\ref{subsect:cech}.
Thus
$\diam f'^{-1}(\si)\le 4d+\diam(Z_\si\times\{t_k\})
 \le\const(\la)$,
and
$f$
is uniformly cobounded. This completes the proof
of Theorem~\ref{thm:conedim}.

\section{Hyperbolic spaces}\label{sect:hypspaces}

\subsection{Basics of hyperbolic spaces}

We briefly recall necessary facts from the hyperbolic spaces
theory. For more details the reader may consult e.g. \cite{BoS}.

Let
$X$
be a metric space. Fix a base point
$o\in X$
and for
$x$, $x'\in X$
put
$(x|x')_o=\frac{1}{2}(|xo|+|x'o|-|xx'|)$.
The number
$(x|x')_o$
is nonnegative by the triangle inequality, and it is
called the Gromov product of
$x$, $x'$
w.r.t.
$o$.
A {\em $\de$-triple} is a triple of three real numbers
$a$, $b$, $c$
with the property that the two smallest of these numbers
differ by at most
$\de$.

A metric space
$X$
is {\em (Gromov) hyperbolic} if for some
$\de\ge 0$,
some base point
$o\in X$
and all
$x$, $x'$, $x''\in X$
the numbers
$(x|x')_o$, $(x'|x'')_o$, $(x|x'')_o$
form a
$\de$-triple.
This condition is equivalent to the
$\de$-{\em inequality}
$$(x|x'')_o\ge\min\{(x|x')_o,(x'|x'')_o\}-\de.$$

Let
$X$
be a hyperbolic space and
$o\in X$
be a base point. A sequence of points
$\{x_i\}\sub X$
{\em converges to infinity,} if
$$\lim_{i,j\to\infty}(x_i|x_j)_o=\infty.$$
Two sequences
$\{x_i\}$, $\{x_i'\}$
that converge to infinity are {\em equivalent} if
$$\lim_{i\to\infty}(x_i|x_i')_o=\infty.$$

The {\em boundary at infinity}
$\di X$
of
$X$
is defined as the set of equivalence classes
of sequences converging to infinity.
The Gromov product extends to
$X\cup\di X$
as follows. For points
$\xi$, $\xi'\in\di X$
the Gromov product is defined by
$$(\xi|\xi')_o=\inf\liminf_{i\to\infty}(x_i|x_i')_o,$$
where the infimum is taken over all sequences
$\{x_i\}\in\xi$, $\{x_i'\}\in\xi'$.
Note that
$(\xi|\xi')_o$
takes values in
$[0,\infty]$,
and that
$(\xi|\xi')_o=\infty$
if and only if
$\xi=\xi'$.
Furthermore, for
$\xi$, $\xi'$, $\xi''\in\di X$
the following holds, see \cite[Sect.3]{BoS}

\begin{itemize}
\item[(1)] for sequences
$\{y_i\}\in \xi$, $\{y'_i\}\in \xi'$
we have
$$(\xi|\xi')_o\leq \liminf_{i\to\infty}{(y_i|y'_i)_o}
\leq \limsup_{i\to\infty}{(y_i|y'_i)_o} \leq (\xi|\xi')_o + 2\de$$

\item[(2)] $(\xi|\xi')_o$, $(\xi'|\xi'')_o$, $(\xi'|\xi'')_o$
 is a
 $\de$-triple.
\end{itemize}

Similarly, the Gromov product
$$(x|\xi)_o=\inf\liminf_{i\to\infty}(x|x_i)_o$$
is defined for any
$x\in X$, $\xi\in\di X$,
where the infimum is taken over all sequences
$\{x_i\}\in\xi$,
and the
$\de$-inequality
holds for any three points from
$X\cup\di X$.

A metric
$d$
on the boundary at infinity
$\di X$
of
$X$
is said to be {\em visual}, if there are
$o\in X$, $a>1$
and positive constants
$c_1$, $c_2$,
such that
$$c_1a^{-(\xi|\xi')_o}\le d(\xi,\xi')\le c_2a^{-(\xi|\xi')_o}$$
for all
$\xi$, $\xi'\in\di X$.
In this case we say that
$d$
is the visual metric w.r.t. the base point
$o$
and the parameter
$a$.

\subsection{The hyperbolic cone}

\begin{pro}\label{pro:hypcone} Let
$Z$
be a bounded metric space. Then the hyperbolic cone
$Y=\cone(Z)$
is a
$\de$-hyperbolic
space with
$\de=\de(\hyp^2)$,
there is a canonical inclusion
$Z\sub\di Y$,
and the metric of
$Z$
is visual. If in addition
$Z$
is complete then
$\di Y=Z$.
\end{pro}

\begin{proof} We can assume that
$(y|y')_o\le(y|y'')_o\le(y''|y')_o$
for
$y$, $y'$, $y''\in Y$.
We show that
$(y|y'')_o\le(y|y')_o+\de$.
To this end, consider triangles
$\ov o\,\ov y\,\ov y''$
and
$\ov o\,\ov y''\,\ov y'$
in
$\hyp^2$
with common side
$\ov o\,\ov y''$
separating them such that
$|\ov o\,\ov y|=|oy|$,
$|\ov o\,\ov y'|=|oy'|$,
$|\ov o\,\ov y''|=|oy''|$,
and
$|\ov y\,\ov y''|=|yy''|$, $|\ov y''\,\ov y'|=|y''y'|$.
Then
$|yy'|\le|\ov y\,\ov y'|$
by the triangle inequality in
$Z$.
It follows that
$(\ov y|\ov y'')_{\ov o}=(y|y'')_o$,
$(\ov y''|\ov y')_{\ov o}=(y''|y')_o$
and
$(\ov y|\ov y')_{\ov o}\le(y|y')_o$.
Therefore,
$(y|y'')_o-(y|y')_o
 \le(\ov y|\ov y'')_{\ov o}-(\ov y|\ov y')_{\ov o}\le\de$
since
$\hyp^2$
is
$\de$-hyperbolic.

For every
$z\in Z$
the ray
$\{z\}\times[0,\infty)\sub Y$
represents a point from
$\di Y$
which we identify with
$z$.
This yields the inclusion
$Z\sub\di Y$.
The last assertion of the proposition is
easy to check.

It remains to show that the metric of
$Z$
is visual. Given
$z$, $z'\in Z$,
consider the geodesic rays
$\ga(t)=(z,t)$, $\ga'(t)=(z',t)$
in
$\cone(Z)$.
Then
$\ga\in z$, $\ga'\in z'$
viewed as points at infinity, and for
$(\ga|\ga')_o=\lim_{t\to\infty}(\ga(t)|\ga'(t))_o$
(it is easy to see that the Gromov product
$(\ga(t)|\ga'(t))_o$
is monotone) we have
$$(z|z')_o\le(\ga|\ga')_o\le(z|z')_o+2\de.$$
For comparison geodesic rays
$\ov\ga$, $\ov\ga'\sub\hyp^2$
with common vertex
$\ov o$
and
$$\angle_{\ov o}(\ov\ga,\ov\ga')=\mu|zz'|$$
(recall
$\mu=\pi/\diam Z$)
we have
$(\ov\ga|\ov\ga')_{\ov o}=(\ga|\ga')_o$
and
$(\ov\ga|\ov\ga')_{\ov o}\le d\le(\ov\ga|\ov\ga')_{\ov o}+\de$,
where
$d=\dist(\ov o,\ov z\,\ov z')$
and
$\ov z\,\ov z'\sub\hyp^2$
is the infinite geodesic with the end points at infinity
$\ov z=\ov\ga(\infty)$, $\ov z'=\ov\ga'(\infty)$.
By the angle of parallelism formula from geometry of
$\hyp^2$
we have
$\tan\frac{\mu|zz'|}{4}=e^{-d}$,
therefore, we conclude that
$$e^{-3\de}e^{-(z|z')_o}\le\tan\frac{\mu|zz'|}{4}\le e^{-(z|z')_o}$$
for every
$z$, $z'\in Z$.
The function
$s\mapsto\frac{1}{s}\tan\frac{\mu s}{4}$
is uniformly bounded and separated from zero
on
$[0,\diam Z]$.
It follows that the metric of
$Z\sub\di Y$
is visual w.r.t. the vertex
$o\in Y$
and the parameter
$a=e$.
\end{proof}

Let
$X$
be a hyperbolic space,
$x_0\in X$
a base point. For
$x\in X$
we denote
$|x|=|xx_0|$.
We also omit the subscript
$x_0$
from the notations of Gromov products w.r.t.
$x_0$.
The space
$X$
is called {\em visual}, if for some base point
$x_0\in X$
there is a positive constant
$D$
such that for every
$x\in X$
there is
$\xi\in\di X$
with
$|x|\le(x|\xi)+D$
(one easily sees that this property is independent of
the choice of
$x_0$).
This definition is due to V.~Schroeder, cf. \cite[Sec.~5]{BoS}.

\begin{pro}\label{pro:roughsim} Every visual hyperbolic space
$X$
is roughly similar to a subspace of the hyperbolic cone
over the boundary at infinity,
$\cone(\di X)$,
where
$\di X$
is taken with a visual metric.
\end{pro}

\begin{proof} We fix a visual metric on
$\di X$
w.r.t.
$x_0\in X$
and a parameter
$a>1$.
Replacing
$X$
by
$\la X$
with
$\la=1/\ln a$
we can assume that
$a=e$.
Since
$X$
is visual, there is a constant
$D>0$
such that for every
$x\in X$
there is
$\xi=\xi(x)\in\di X$
with
$|x|\le(x|\xi)+D$.
We define
$F:X\to Y$, $Y=\cone(\di X)$
by
$F(x)=(\xi(x),|x|)\in Y$.
Note that
$F(x_0)=o$.

It follows from Proposition~\ref{pro:hypcone} that
the Gromov product
$(\xi|\xi')$
in
$X$
coincides with the Gromov product
$(\xi|\xi')_o$
in
$Y$
up to a uniformly bounded error for every
$\xi$, $\xi'\in\di X$.
Since
$|F(x)|=|x|$
for every
$x\in X$,
by \cite[Lemma~5.1]{BoS} we have
$$|xx'|\doteq|x|+|x'|-2\min\{(\xi(x)|\xi(x')),|x|,|x'|\}
  \doteq|F(x)F(x')|$$
up uniformly bounded error for every
$x$, $x'\in X$.
Hence
$F$
is roughly isometric and
$X$
is roughly similar to a subspace of
$Y$.
\end{proof}

\subsection{Proof of Theorem~\ref{thm:main}} The asymptotic
dimension is a quasi-isometry invariant, thus
$\asdim X\le\asdim\cone(\di X)$
by Proposition~\ref{pro:roughsim}. By Theorem~\ref{thm:conedim}
we obtain
$\asdim X\le\cdim(\di X)+1$.
\qed

\begin{rem}\label{rem:gromov} Coming back to the Gromov
argument (see Introduction) that
$\asdim X\le n$
for every negatively pinched Hadamard manifold
$X$
of dimension
$n$,
in my opinion to complete the proof one needs to show
that
$\cdim\di X=n-1$
where
$\di X=S^{n-1}$
is considered with a visual metric. However, at the moment
this is only known for
$X=\hyp^n$.
\end{rem}

%%%%%%%%%%%%%%%%%%%%%%%%%%%%%%%%%%%%%%%%%%%%%%%%%%%%%%%%%%%%%%

%\bigskip
%\begin{tabbing}
%
%
%
%St. Petersburg Dept. of Steklov\hskip11em\relax \= \\
%
%Math. Institute RAS, Fontanka 27, \> \\
%
%191023 St. Petersburg, Russia\> \\
%
%{\tt sbuyalo@pdmi.ras.ru}\> \\
%
%\end{tabbing}

\end{document}